\documentclass{amsart}
\usepackage{amssymb}
\theoremstyle{plain}
 \newtheorem{thm}{Theorem}
 
 \newtheorem{lem}{Lemma}
 
\theoremstyle{definition}

 \newtheorem{exmp}{Example}
\theoremstyle{remark}
 \newtheorem{rem}{Remark\ignorespaces}
 
\newcommand{\NaturalNumber}{\mathbb N}
\newcommand{\RealNumber}{\mathbb R}

\renewcommand{\labelenumi}{(\roman{enumi})}
\allowdisplaybreaks
\begin{document}
\title[Example for nonexpansive semigroup]
{An example for a one-parameter nonexpansive semigroup}
\author[T. Suzuki]{Tomonari Suzuki}
\date{}
\hyphenation{kyu-shu kita-kyu-shu to-bata-ku sen-sui-cho}
\address{
Department of Mathematics,
Kyushu Institute of Technology,
1-1, Sensuicho, Tobataku, Kitakyushu 804-8550, Japan}
\email{suzuki-t@mns.kyutech.ac.jp}

\begin{abstract}
In this paper,
 we give one example for a one-parameter nonexpansive semigroup.
This Example shows that
 there exists a one-parameter nonexpansive semigroup
 $\{ T(t) : t \geq 0 \}$ on a closed convex subset $C$ of a Banach space $E$
 such that
 $$ \lim_{t \rightarrow \infty}
 \left\| \frac{1}{t} \int_0^t T(s)x \; ds - x \right\| = 0 $$
 for some $x \in C$,
 which is not a common fixed point of $\{ T(t) : t \geq 0 \}$.
\end{abstract}
\maketitle

\section{Introduction}
\label{SC:introduction}

Throughout this paper,
 we denote by $\NaturalNumber$ and $\RealNumber$
 the set of positive integers and real numbers, respectively.

A family $\{ T(t): t \geq 0 \}$ of mappings on $C$
 is called a one-parameter nonexpansive semigroup
 on a subset $C$ of a Banach space $E$ if the following hold:
\begin{enumerate}
\renewcommand{\labelenumi}{(sg\arabic{enumi})}
\renewcommand{\theenumi}{(sg\arabic{enumi})}
\item\label{ENUM:sg:nonex}
 For each $t \geq 0$, $T(t)$ is a nonexpansive mapping on $C$, i.e.,
 $$ \| T(t)x - T(t)y \| \leq \| x - y \| $$
 for all $x,y \in C$;
\item\label{ENUM:sg:T0}
 $T(0)x = x$ for all $x \in C$;
\item\label{ENUM:sg:s+t}
 $T(s+t) = T(s) \circ T(t)$ for all $s, t \geq 0$;
\item\label{ENUM:sg:conti}
 For each $x \in C$, the mapping $t \mapsto T(t)x$ is continuous.
\end{enumerate}
We know that
 $\{ T(t) : t \geq 0 \}$ has a common fixed point
 under the assumption that
 $C$ is weakly compact convex and $E$ has the Opial property;
 see \cite{REF:Belluce_Kirk1967_Illinois,
 REF:Browder1965_ProcNAS_3,
 REF:Bruck1974_Pacific,
 REF:DeMarr1963_Pacific,
 REF:Gossez_LamiDazo1972_Pacific,
 REF:Lim1974_Pacific,
 REF:Opial1967_BullAMS} and others.

Convergence theorems for one-parameter nonexpansive semigroups
 are proved in
 \cite{REF:Baillon1976,
 REF:Baillon_Brezis1976_Houston,
 REF:Hirano1982_JMSJapan,
 REF:Miyadera_Kobayasi1982_NATMA,
 REF:TS2003_ProcAMS,
 REF:TSP_mcbt_1_04}
 and others.
For example,
 Baillon and Brezis in \cite{REF:Baillon_Brezis1976_Houston}
 proved the following;
 see also page 80 in \cite{REF:Takahashi_ybook}.

\begin{thm}[Baillon and Brezis \cite{REF:Baillon_Brezis1976_Houston}]
\label{THM:Baillon-Brezis}
Let $C$ be a bounded closed convex subset of a Hilbert space $E$ and
 let $\{ T(t) : t \geq 0 \}$ be
 a one-parameter nonexpansive semigroup on $C$.
Then, for any $x \in C$,
 $$ \frac{1}{t} \int_{0}^{t} T(s) x \; ds $$
 converges weakly to a common fixed point of $\{ T(t) : t \geq 0 \}$
 as $t \rightarrow \infty$.
\end{thm}

\noindent
Also, Suzuki and Takahashi in \cite{REF:TSP_mcbt_1_04}
 proved the following.

\begin{thm}[Suzuki and Takahashi \cite{REF:TSP_mcbt_1_04}]
\label{THM:TS-Takahashi-conv}
Let $C$ be a compact convex subset of a Banach space $E$ and
 let $\{ T(t) : t \geq 0 \}$ be
 a one-parameter nonexpansive semigroup on $C$.
Let $x_1 \in C$ and
 define a sequence $\{ x_n \}$ in $C$ by
 $$ x_{n+1} =  \frac{\alpha_n}{t_n} \int_{0}^{t_n} T(s) x_n \; d s
  +(1 - \alpha_n) x_n $$
 for $n \in \NaturalNumber$,
 where $\{ \alpha_n \} \subset [0,1]$ and $\{ t_n \} \subset (0,\infty)$
 satisfy the following conditions:
 $$ 0 < \liminf_{n \rightarrow \infty} \alpha_n \leq
 \limsup_{n \rightarrow \infty} \alpha_n < 1, \quad
 \lim_{n \rightarrow \infty} t_n = \infty,
 \quad\text{and}\quad
 \lim_{n \rightarrow \infty} \frac{t_{n+1}}{t_n} = 1 . $$
Then $\{ x_n \}$ converges strongly to
 a common fixed point $z_0$ of $\{ T(t) : t \geq 0 \}$.
\end{thm}

\noindent
In the proof of it, the following Theorem plays a very important role.

\begin{thm}[Suzuki and Takahashi \cite{REF:TSP_mcbt_1_04}]
\label{THM:int-compact}
Let $C$ be a compact convex subset of a Banach space $E$.
Let $\{ T(t): t \geq 0 \}$ be a one-parameter nonexpansive semigroup on $C$.
Then for $z \in C$,
 the following are equivalent:
\begin{enumerate}
\item
 $z$ is a common fixed point of $\{ T(t) : t \geq 0 \}$;
\item
 $$ \liminf_{t \rightarrow \infty} \left\|
  \frac{1}{t} \int_{0}^{t} T(s) z \; d s - z
 \right\| = 0 $$
 holds.
\end{enumerate}
\end{thm}

\noindent
Recently, Suzuki proved in \cite{REF:Suzuki-thm}
 the following similar result to Theorem \ref{THM:int-compact}.
This Theorem also plays a very important role
 in the proof of the existence of nonexpansive retraction
 onto the set of common fixed points.

\begin{thm}[Suzuki \cite{REF:Suzuki-thm}]
\label{THM:int}
Let $E$ be a Banach space with the Opial property and
 let $C$ be a weakly compact convex subset of $E$.
Let $\{ T(t): t \geq 0 \}$ be a one-parameter nonexpansive semigroup on $C$.
Then for $z \in C$,
 the following are equivalent:
\begin{enumerate}
\item
 $z$ is a common fixed point of $\{ T(t) : t \geq 0 \}$;
\item
 $$ \liminf_{t \rightarrow \infty} \left\|
  \frac{1}{t} \int_{0}^{t} T(s) z \; d s - z
 \right\| = 0 $$
 holds;
\item
 there exists a subnet of
 a net
 $$ \left\{ \frac{1}{t} \int_{0}^{t} T(s) z \; ds  \right\} $$
 in $C$
 converging weakly to $z$.
\end{enumerate}
\end{thm}

So, it is natural problem
 whether or not
 the conclusion of Theorems \ref{THM:int-compact} and \ref{THM:int}
 holds in general.
In this paper,
 we give one example concerning Theorems \ref{THM:int-compact}
 and \ref{THM:int}.
This Example shows that
 there exists a one-parameter nonexpansive semigroup
 $\{ T(t) : t \geq 0 \}$ on a closed convex subset $C$ of a Banach space $E$
 such that
 $$ \lim_{t \rightarrow \infty}
 \left\| \frac{1}{t} \int_0^t T(s)x \; ds - x \right\| = 0 $$
 for some $x \in C$,
 which is not a common fixed point of $\{ T(t) : t \geq 0 \}$.
That is, our answer of the problem is negative.

\section{Example}
\label{SC:example}

We give one Example concerning Theorems \ref{THM:int-compact} and
 \ref{THM:int}.
See also Example 3.7 in \cite{REF:Goebel_Kirk}.

\begin{exmp}
\label{EX:int}
Put $\Omega = \{ -1 \} \cup [0,\infty)$,
 let $E$ be the Banach space consisting of all bounded continuous functions
 on $\Omega$ with supremum norm, and
 define a subset $C$ of $E$ by
 $$ C = \left\{ x \in E :
  \begin{array}{l}
   0 \leq x(u) \leq 1 \quad \text{for } u \in \Omega \\
   | x(u_1) - x(u_2) | \leq | u_1 - u_2 | \quad
    \text{for } u_1, u_2 \in [0,\infty)
  \end{array}
 \right\} . $$
Define a nonexpansive semigroup $\{ T(t) : t \geq 0 \}$ as follows:
For $t \in [0,1]$, define
 $$ \big( T(t)x \big) (u)
 =
 \begin{cases}
  x(u), & \text{if $u = -1$}, \\
  x(u-t), & \text{if $u \geq t$}, \\
  x(0) - t + u, & \text{if $0 \leq u \leq t$}, \\
  & 1 - \alpha_x(1-t+u) \leq x(0)-t+u, \\
  x(0) + t - u, & \text{if $0 \leq u \leq t$}, \\
  & 1 - \alpha_x(1-t+u) \geq x(0)+t-u, \\
  1 - \alpha_x(1-t+u), & \text{if  $0 \leq u \leq t$}, \\
  & \big| 1 - \alpha_x(1-t+u) - x(0) \big| \leq t - u,
 \end{cases} $$
 where
 $$ \alpha_x(1-t+u) =
  \sup\big\{ x(s) : s \in \{ -1 \} \cup [1-t+u,\infty) \big\} .$$
For $t \in (1,\infty)$,
 there exist $m \in \NaturalNumber$ and
 $t' \in [0,1/2)$ satisfying $ t = m / 2 + t'$.
Define $T(t)$ by
 $$ T(t) = T(1/2)^m \circ T(t') .$$
Then
 $0 \in C$ is not a common fixed point of $\{ T(t) : t \geq 0 \}$ and
\begin{equation}
\label{EQU:int:int}
 \lim_{t \rightarrow \infty} \left\|
 \frac{1}{t} \int_{0}^{t} T(s) 0 \; d s - 0 \right\| = 0
\end{equation}
 holds.
\end{exmp}

Before proving Example \ref{EX:int},
 we need some lemmas.

\begin{lem}
\label{LEM:int:alpha}
 The following hold:
 \begin{enumerate}
 \item
  $ | \alpha_x(u_1) - \alpha_x(u_2) | \leq | u_1 - u_2 | $
  for $x \in C$ and $u_1, u_2 \in [0,\infty)$;
 \item
  $ | \alpha_x(u) - \alpha_y(u) | \leq \| x - y \| $
  for $x,y \in C$ and $u \in [0,\infty)$.
 \end{enumerate}
\end{lem}

\begin{proof}
We first show (i).
Without loss of generality, we may assume $u_1 < u_2$.
For $s \in [u_1,u_2]$,
 we have $| x(s) - x(u_2) | \leq | s - u_2 |$ and hence
 $$ x(s)
 \leq x(u_2) + | s - u_2 |
 \leq \alpha_x(u_2) + | u_1 - u_2 | . $$
For $s \in [u_2,\infty)$, we have
 $$ x(s) \leq \alpha_x(u_2) \leq \alpha_x(u_2) + | u_1 - u_2 | . $$
Hence,
 $$ \alpha_x(u_1)
 \leq \alpha_x(u_2) + | u_1 - u_2 | $$
 holds.
Since $ \alpha_x(u_2) \leq \alpha_x(u_1) $,
 we obtain
 $$ | \alpha_x(u_1) - \alpha_x(u_2) | \leq | u_1 - u_2 | .$$
We next show (ii).
For each $\varepsilon > 0$,
 there exists $s \in \{ -1 \} \cup [ u, \infty)$ satisfying
 $ x(s) > \alpha_x(u) - \varepsilon $.
We have
 $$ \alpha_x(u) - \alpha_y(u)
 \leq x(s) + \varepsilon - y(s)
 \leq \| x - y \| + \varepsilon .$$
Since $\varepsilon$ is arbitrary,
 we have
 $ \alpha_x(u) - \alpha_y(u) \leq \| x - y \| $.
Similarly we obtain
 $ \alpha_y(u) - \alpha_x(u) \leq \| x - y \| $
 and hence
 $ | \alpha_x(u) - \alpha_y(u) | \leq \| x - y \| $.
\end{proof}

\begin{lem}
\label{LEM:int:3}
Fix $x \in C$, $t \in [0,1]$, and $u_1, u_2$ with $0 \leq u_1 \leq u_2 \leq t$.
Then the following hold:
 \begin{enumerate}
 \item
 $1 - \alpha_x(1-t+u_1) < \big( T(t)x \big)(u_2) - u_2 + u_1$
 implies
 $\big( T(t)x \big)(u_1) = x(0)-t+u_1$ and
 $\big( T(t)x \big)(u_2) = x(0)-t+u_2$;
 \item
 $1 - \alpha_x(1-t+u_1) > \big( T(t)x \big)(u_2) + u_2 - u_1$
 implies
 $\big( T(t)x \big)(u_1) = x(0)+t-u_1$ and
 $\big( T(t)x \big)(u_2) = x(0)+t-u_2$;
 \item
 $\left| 1 - \alpha_x(1-t+u_1) - \big( T(t)x \big)(u_2) \right|
 \leq u_2 - u_1$
 implies
 $\big( T(t)x \big)(u_1) = 1-\alpha_x(1-t+u_1)$.
 \end{enumerate}
\end{lem}

\begin{rem}
One and only one of the assumptions of (i), (ii) and (iii) holds.
\end{rem}

\begin{proof}
We first prove (i).
We assume that
 $ 1 - \alpha_x(1-t+u_2) > x(0) - t + u_2$.
Then by the definition of $T(t)$,
 $$ \big( T(t)x \big)(u_2)
 = \min\{ x(0)+t-u_2, 1-\alpha_x(1-t+u_2) \} .$$
So, we have
 \begin{align*}
 \big( T(t)x \big)(u_2) - u_2 + u_1
 &\leq 1-\alpha_x(1-t+u_2) - u_2 + u_1 \\*
 &\leq 1-\alpha_x(1-t+u_1)
 \end{align*}
 by Lemma \ref{LEM:int:alpha}.
This is a contradiction.
Therefore we obtain
 $ 1 - \alpha_x(1-t+u_2) \leq x(0) - t + u_2$.
Hence $\big( T(t)x \big)(u_2) = x(0)-t+u_2$.
Since
 $$ 1 - \alpha_x(1-t+u_1) < \big( T(t)x \big)(u_2) - u_2 + u_1
 = x(0) - t + u_1 ,$$
 we have $\big( T(t)x \big)(u_1) = x(0)-t+u_1$.
Similarly, we can prove (ii).
We finally prove (iii).
We assume that
 $ 1 - \alpha_x(1-t+u_1) < x(0) - t + u_1$.
Then by Lemma \ref{LEM:int:alpha}, we have
 \begin{align*}
 1 - \alpha_x(1-t+u_2)
 &\leq 1 - \alpha_x(1-t+u_1) + u_2 - u_1 \\*
 &< x(0) - t + u_1 + u_2 - u_1
 = x(0) - t + u_2 .
 \end{align*}
Hence $\big( T(t)x \big)(u_2) = x(0) - t + u_2$.
So,
 \begin{align*}
 \big( T(t)x \big)(u_2) - \big( 1 - \alpha_x(1 - t + u_1) \big)
 &> \big( x(0) - t + u_2 \big) - \big( x(0) - t + u_1 \big) \\*
 &= u_2 - u_1 .
 \end{align*}
This is a contradiction.
Therefore we obtain
 $ 1 - \alpha_x(1-t+u_1) \geq x(0) - t + u_1$.
Similarly we can prove
 $ 1 - \alpha_x(1-t+u_1) \leq x(0) + t - u_1$.
Hence $\big( T(t)x \big)(u_1) = 1-\alpha_x(1-t+u_1)$.
\end{proof}

\begin{proof}[Proof of Example \ref{EX:int}]
It is clear that $C$ is closed and convex.
We first prove that
 $T(t) x \in C$ for all $t \in [0,1]$ and $x \in C$.
It is clear that
 $$ 0 \leq \big( T(t) x \big)(-1) = x(-1) \leq 1 $$
 and
 $$ 0 \leq \big( T(t) x \big)(u) = x(u-t) \leq 1 $$
 for $u \in [t,\infty)$.
For $u \in [0,t]$, since
 $ 0 \leq 1 - \alpha_x(1-t+u) \leq 1 $,
 $ x(0)-t+u \leq x(0) \leq 1$ and
 $ x(0)+t-u \geq x(0) \geq 0$,
 we have
 $ 0 \leq \big( T(t) x \big)(u) \leq 1 $.
Fix $u_1, u_2 \in [0,\infty)$ with $u_1 < u_2$.
In the case of $t \leq u_1$,
 we have
\begin{align*}
 \left| \big( T(t)x \big)(u_1) - \big( T(t)x \big)(u_2) \right|
 &= | x(u_1 - t) - x(u_2 - t) | \\*
 &\leq | (u_1 - t) - (u_2 - t) |
 = | u_1 - u_2 | .
 \end{align*}
In the case of $u_2 \leq t$, by Lemma \ref{LEM:int:3},
 it is easily proved that
 $$ \left| \big( T(t)x \big)(u_1) - \big( T(t)x \big)(u_2) \right|
 \leq | u_1 - u_2 | .$$
In the case of $u_1 \leq t \leq u_2$,
 we have
\begin{align*}
 & \left| \big( T(t)x \big)(u_1) - \big( T(t)x \big)(u_2) \right| \\*
 &\leq \left| \big( T(t)x \big)(u_1) - \big( T(t)x \big)(t) \right|
  + \left| \big( T(t)x \big)(t) - \big( T(t)x \big)(u_2) \right| \\*
 &\leq | u_1 - t | + | t - u_2 |
 = | u_1 - u_2 | .
 \end{align*}
Therefore we have shown $T(t)x \in C$ for $t \in [0,1]$ and $x \in C$.
By the definition of $\{ T(t) : t \geq 0 \}$,
 we have
 $T(t) x \in C$ for all $t \in [0,\infty)$ and $x \in C$.
We next show that
 $\{ T(t) : t \geq 0 \}$ is a one-parameter nonexpansive semigroup on $C$.

\ref{ENUM:sg:nonex}:
Fix $t \in [0,1]$, and $x, y \in C$.
We shall prove
 \begin{equation}
 \label{EQU:int:nonex}
 \Big| \big( T(t) x \big) (u) - \big( T(t) y \big) (u) \Big| \leq
 \| x - y \|
 \end{equation}
 for all $u \in \Omega$.
We have
 $$ \big| \big( T(t) x \big) (-1) - \big( T(t) y \big)(-1) \big|
 = \big| x(-1) - y(-1) \big|
 \leq \| x - y \| .$$
For $u \geq t$, we have
 $$ \big| \big( T(t)x \big)(u) - \big( T(t)y \big)(u) \big|
 = \big| x(u-t) - y(u-t) \big|
 \leq \| x - y \| .$$
Fix $u$ with $0 \leq u \leq t$.
In the case of
 $ 1 - \alpha_x(1-t+u) \leq x(0)-t+u $ and
 $ 1 - \alpha_y(1-t+u) \leq y(0)-t+u $,
 we have
 \begin{align*}
 \left| \big( T(t) x \big) (u) - \big( T(t) y \big) (u) \right|
 &= \left| \big( x(0) - t + u \big) - \big( y(0) - t + u \big) \right| \\*
 &= \left| x(0) - y(0) \right|
 \leq \| x - y \| .
 \end{align*}
In the case of
 $ 1 - \alpha_x(1-t+u) \leq x(0)-t+u $ and
 $ 1 - \alpha_y(1-t+u) > y(0)-t+u $,
 we have
 $$ \big( T(t) y \big)(u) = \min\big\{ 1 - \alpha_y(1-t+u), y(0)+t-u \big\}
 \geq y(0)-t+u . $$
Hence,
 \begin{align*}
 \big( T(t) x \big) (u) - \big( T(t) y \big) (u)
 &\leq \big( x(0) - t + u \big) - \big( y(0) - t + u \big) \\*
 &= x(0) - y(0)
 \leq \| x - y \|
 \end{align*}
 and
 \begin{align*}
 \big( T(t) y \big) (u) - \big( T(t) x \big) (u)
 &\leq \big( 1 - \alpha_y(1-t+u) \big) - \big( 1 - \alpha_x(1-t+u) \big) \\*
 &= \alpha_x(1-t+u) - \alpha_y(1-t+u)
 \leq \| x - y \|
 \end{align*}
 hold.
Therefore \eqref{EQU:int:nonex} holds.
Similarly we can prove \eqref{EQU:int:nonex} in the other cases.
On the other hand,
 we have
 \begin{align*}
 \| T(t)x - T(t)y \|
 &\geq \sup\left\{ \left| \big( T(t)x \big)(u) - \big( T(t)y \big)(u) \right|
  : u \in \{-1\} \cup [t,\infty) \right\} \\*
 &= \sup\{ | x(u) - y(u) : u \in \Omega \}
 = \| x - y \| .
\end{align*}
Hence we have shown
 \begin{equation}
 \label{EQU:int:isometric}
 \left\| T(t) x - T(t)y \right\| = \| x - y \|
 \end{equation}
 for $t \in [0,1]$ and $x,y \in C$.
So, by the definition of $\{ T(t) : t \geq 0 \}$,
 \eqref{EQU:int:isometric} hold for all $t \in [0,\infty)$
 and $x,y \in C$.

\ref{ENUM:sg:T0}:
It is clear that
 $T(0)$ is the identity mapping on $C$.

\ref{ENUM:sg:s+t}:
Fix $t_1, t_2 \in [0,1/2]$ and $x \in C$.
We shall prove that
 \begin{equation}
 \label{EQU:int:s+t}
 \big( T(t_1) \circ T(t_2) x \big) (u) = \big( T(t_1 + t_2) x \big)(u)
 \end{equation}
 for all $u \in \Omega$.
We have
 $$ \big( T(t_1) \circ T(t_2) x \big) (-1)
 =  \big( T(t_2) x \big) (-1)
 = x(-1)
 = \big( T(t_1 + t_2) x \big)(-1) . $$
For $u \geq t_2$, we have
 \begin{align*}
 \big( T(t_1+t_2)x \big)(t_1 + u)
 &= x \big( (t_1+u) - (t_1+t_2) \big)
 = x(u - t_2) \\*
 &= \big( T(t_2)x \big)(u) .
 \end{align*}
For $u \in [0,t_2]$,
 since
  $t_1 + u \leq t_1 + t_2$,
  $1 - \alpha_x(1-t_2+u) = 1 - \alpha_x \big( 1 - (t_1+t_2) + (t_1+u) \big)$,
  $x(0) - t_2 + u = x(0) - (t_1+t_2) + (t_1+u)$ and
  $x(0) + t_2 - u = x(0) + (t_1+t_2) - (t_1+u)$,
 the two definitions of $\big( T(t_1+t_2)x \big)(t_1+u)$ and
 $\big( T(t_2)x \big)(u)$ coincide.
Therefore
 $$ \big( T(t_1+t_2)x \big)(t_1+u) = \big( T(t_2)x \big)(u) . $$
So, for $u \geq t_1$,
 \begin{align*}
 \big( T(t_1) \circ T(t_2) x \big)(u)
 &= \big( T(t_2)x \big)(u-t_1)
 = \big( T(t_1+t_2) x \big) \big(t_1 + (u-t_1) \big) \\*
 &= \big( T(t_1+t_2)x \big) (u) .
 \end{align*}
Fix $u$ with $0 \leq u \leq t_1$.
Then we have
 \begin{align*}
 & 1 - \alpha_{T(t_2)x}(1-t_1+u) \\*
 &= 1 - \sup\left\{ \big( T(t_2)x \big)(s)
  : s \in \{ -1 \} \cup [1-t_1+u,\infty)\right\} \\
 &= 1 - \max \Big\{ x(-1), \sup\left\{ x(s-t_2)
  : s \in [1-t_1+u,\infty)\right\} \Big\} \\*
 &= 1 - \alpha_x(1-t_1-t_2+u) .
 \end{align*}
In the case of
 $1 - \alpha_{T(t_2)x}(1-t_1+u) < \big( T(t_2)x \big)(0) - t_1 + u$,
 we have
 $$ \big( T(t_1) \circ T(t_2)x \big)(u)
 = \big( T(t_2)x \big)(0) - t_1 + u . $$
Since
 \begin{align*}
 1 - \alpha_x(1-t_1-t_2+u)
 &= 1 - \alpha_{T(t_2)x}(1-t_1+u) \\*
 &< \big( T(t_2)x \big)(0) - t_1 + u \\
 &= \big( T(t_1+t_2)x \big)(t_1) - t_1 + u,
 \end{align*}
 we have
 $$ \big( T(t_1+t_2)x \big)(u) = x(0) - t_1 - t_2 + u $$
 and
 $$ \big( T(t_1+t_2)x \big)(t_1) = x(0) - t_1 - t_2 + t_1
 = x(0) - t_2 $$
 by Lemma \ref{LEM:int:3}.
So,
 \begin{align*}
 \big( T(t_1) \circ T(t_2)x \big)(u)
 &= \big( T(t_2)x \big)(0) - t_1 + u \\*
 &= \big( T(t_1+t_2)x \big)(t_1) - t_1 + u \\
 &= x(0) - t_2 - t_1 + u \\*
 &= \big( T(t_1+t_2)x \big)(u) .
 \end{align*}
Similarly, we can prove
 $ \big( T(t_1) \circ T(t_2)x \big)(u) = \big( T(t_1+t_2)x \big)(u) $
 in the cases of
 $1 - \alpha_{T(t_2)x}(1-t_1+u) > \big( T(t_2)x \big)(0) + t_1 - u$ and
 $ \big| 1 - \alpha_{T(t_2)x}(1-t_1+u) - \big( T(t_2)x \big)(0) \big|
 \leq t_1 - u$.
Therefore $T(t_1) \circ T(t_2) = T(t_1 + t_2)$.
So, we have for $t \in [1/2,1)$,
 $$ T(t) = T(1/2) \circ T(t-1/2) \quad\text{and}\quad
 T(1) = T(1/2) \circ T(1/2) \circ T(0) .$$
Fix $t_1, t_2 \in [0,\infty)$.
Then there exist $m_1, m_2 \in \NaturalNumber \cup \{ 0 \}$ and
 $t_1', t_2' \in [0,1/2)$ satisfying
 $t_1 = m_1 / 2 + t_1'$ and $t_2 = m_2 / 2 + t_2'$.
We have
 \begin{align*}
 T(t_1) \circ T(t_2)
 &= T(1/2)^{m_1} \circ T(t_1') \circ T(1/2)^{m_2} \circ T(t_2') \\*
 &= T(1/2)^{m_1+m_2} \circ T(t_1') \circ T(t_2') \\
 &= T(1/2)^{m_1+m_2} \circ T(\min\{ t_1'+t_2', 1/2 \}) \circ T(\max\{ 0, t_1'+t_2'-1/2\}) \\*
 &= T(t_1 + t_2) .
\end{align*}

\ref{ENUM:sg:conti}:
For $x \in C$ and $t \in [0,\infty)$,
 we have
 \begin{align*}
 \| T(t)x - x \|
 &= \sup \left\{ \left| \big( T(t)x \big)(u) - x(u) \right|
  : u \in [0,\infty) \right\} \\*
 &= \sup \left\{ \left| \big( T(t)x \big)(u) - \big( T(t)x \big)(t+u) \right|
  : u \in [0,\infty) \right\} \\*
 &\leq \sup \left\{ \left| u - (t+u) \right|
  : u \in [0,\infty) \right\}
 = t .
 \end{align*}
Therefore
 we obtain
 $$ \| T(t_1) x - T(t_2) x \|
 = \| T(|t_1-t_2|) x - x \|
 \leq | t_1 - t_2 | $$
 for $x \in C$ and $t_1, t_2 \in [0,1]$.
Therefore $T(\cdot)x$ is continuous for all $x \in C$.

Let us prove
\begin{equation}
\label{EQU:int:F}
 \bigcap_{t \geq 0} F \big( T(t) \big) =
 \big\{ v_s : s \in [0, 1/2] \big\} \cup \big\{ w_s : s \in [0, 1/2] \big\} ,
\end{equation}
 where
 $$ v_s(u) =
 \begin{cases}
  1-s, & \text{if $u = -1$}, \\
  s, & \text{if $u \in [0,\infty)$}
 \end{cases} $$
 and
 $$ w_s(u) =
 \begin{cases}
  s, & \text{if $u = -1$}, \\
  1/2, & \text{if $u \in [0,\infty)$}.
 \end{cases} $$
Fix $s \in [0,1/2]$ and $t \in [0,1]$.
Then we have
 $$ \big| 1 - \alpha_{v_s} (1-t+u) - v_s(0) \big|
 = \big| 1 - (1-s) - s \big| = 0 \leq t - u $$
 and
 $$ \big| 1 - \alpha_{w_s} (1-t+u) - w_s(0) \big|
 = \big| 1 - 1/2 - 1/2 \big| = 0 \leq t - u $$
 for $u \in [0,t]$.
So
 $$ \big( T(t) v_s \big) (u) = 1 - \alpha_{v_s}(1-t+u) = s = v_s(u) $$
 and
 $$ \big( T(t) w_s \big) (u) = 1 - \alpha_{w_s}(1-t+u) = 1/2 = w_s(u) .$$
Hence
 $ T(t) v_s = v_s $ and $ T(t) w_s = w_s$.
Therefore
 $ v_s $ and $w_s$ are common fixed points of $\{ T(t) : t \geq 0 \}$.
Conversely, we assume that $x \in C$ is a common fixed point
 of $\{ T(t) : t \geq 0 \}$.
Put $s = x(0)$.
Then we have
 $$ x(t+u) = \big( T(t) x \big) (t+u) = x(t+u-t) = x(u) $$
 for all $u \in [0,\infty)$ and $t \in [0,1]$.
So, $ x(u) = x(0) = s $ hold
 for all $u \in [0,\infty)$.
We also have
 \begin{align*}
 s &= x(0) = \big( T(1) x \big) (0) = 1 - \alpha_x(1-1+0) = 1 - \alpha_x(0) \\*
 &= \min\{ 1 - x(-1), 1 - s \}
 \end{align*}
Hence $x(-1) \leq 1 - s$ and $s \leq 1/2$.
If $s = 1/2$, then $x = w_{x(-1)}$.
If $s < 1/2$, then $x(-1) = 1-s$ and hence $x = v_s$.
Therefore we have shown \eqref{EQU:int:F}.

Define a function $f$ from $\RealNumber$ into $[0,1]$ by
 $$ f(u) =
 \begin{cases}
  0, & \text{if $u \geq 0$}, \\
  - u, & \text{if $-1 \leq u \leq 0$}, \\
  u + 2, & \text{if $-2 \leq u \leq -1$}, \\
  0, & \text{if $u \leq -2$}.
 \end{cases} $$
We finally show
\begin{equation}
\label{EQU:int:T(t)0}
 \big( T(t) 0 \big)(u) =
 \begin{cases}
  0, & \text{if $u = -1$}, \\
  f(u-t), & \text{if $u \in [0,\infty)$}.
 \end{cases}
\end{equation}
Fix $t \in [0,1]$ and $u \in [0,t]$.
Then we have
 $$ 1 - \alpha_0(1-t+u) = 1 \geq 0 + t - u $$
 and hence
 $$\big( T(t) 0 \big) (u) = 0 + t - u = t-u = f(u-t) $$
 because of $-1 \leq u - t \leq 0$.
Therefore
 $$
 \big( T(1) 0 \big)(s) =
 \begin{cases}
  0, & \text{if $s = -1$ or $s \geq 1$}, \\
  1-s, & \text{if $0 \leq s \leq 1$}.
 \end{cases} $$
Since
 \begin{align*}
 1 - \alpha_{T(1)0} (1-t+u)
 &= 1 - (1-(1-t+u))
 = 1-t+u \\
 &= \big( T(1) 0 \big) (0) - t + u,
 \end{align*}
 we have
 \begin{align*}
 \big( T(t+1) 0 \big) (u)
 &= \big( T(t) \circ T(1) 0 \big) (u)
 = \big( T(1) 0 \big) (0) - t + u \\
 &= 1 - t + u = f \big( u - (1+t) \big) .
 \end{align*}
Therefore
 $$
 \big( T(2) 0 \big)(s) =
 \begin{cases}
  0, & \text{if $s = -1$ or $s \geq 2$}, \\
  2-s, & \text{if $1 \leq s \leq 2$}, \\
  s, & \text{if $0 \leq s \leq 1$}.
 \end{cases} $$
Since
 $$
 \Big| 1 - \alpha_{T(2)0} (1-t+u) - \big( T(2)0 \big) (0) \Big|
 = | 1 - 1 - 0 |
 = 0
 \leq t - u,
 $$
 we have
 \begin{align*}
 \big( T(t+2) 0 \big) (u)
 &= \big( T(t) \circ T(2) 0 \big) (u)
 = 1 - \alpha_{T(2)0} (1 - t + u) \\
 &= 0 = f \big( u - (2+t) \big) .
 \end{align*}
Similarly,
 for $k \in \NaturalNumber$ with $k > 2$,
 we can prove
 $$ \big( T(t+k) 0 \big)(u) = 0 = f \big( u - (k+t) \big) .$$
Therefore we have shown \eqref{EQU:int:T(t)0}.
So, \eqref{EQU:int:int} clearly holds.
This completes the proof.
\end{proof}


\begin{thebibliography}{99}

\bibitem{REF:Baillon1976}
 J. B. Baillon,
 {\it ``Quelques properi\'et\`es de convergence asymptotique pour les
 semigroupes de contractions impa{\`\i}res''},
 C. R. Acad.\ Sci.\ Paris, {\bf 283} (1976), 75--78.

\bibitem{REF:Baillon_Brezis1976_Houston}
 J. B. Baillon and H. Br\'{e}zis,
 {\it ``Une remarque sur le comportement asymptotique des semigroupes non
 lineaires''},
 Houston J.\ Math., {\bf 2} (1976), 5--7.

\bibitem{REF:Belluce_Kirk1967_Illinois}
 L. P. Belluce and W. A. Kirk,
 {\it ``Nonexpansive mappings and fixed-points in Banach spaces''},
 Illinois J.\ Math, {\bf 11} (1967), 474--479.

\bibitem{REF:Browder1965_ProcNAS_3}
 F. E. Browder,
 {\it ``Nonexpansive nonlinear operators in a Banach space''},
 Proc.\ Nat.\ Acad.\ Sci.\ USA, {\bf 54} (1965), 1041--1044.

\bibitem{REF:Bruck1974_Pacific}
 R. E. Bruck,
 {\it ``A common fixed point theorem for a commuting family of nonexpansive
 mappings''},
 Pacific J.\ Math., {\bf 53} (1974), 59--71.

\bibitem{REF:DeMarr1963_Pacific}
 R. DeMarr,
 {\it ``Common fixed points for commuting contraction mappings''},
 Pacific J.\ Math., {\bf 13} (1963), 1139--1141.

\bibitem{REF:Goebel_Kirk}
 K. Goebel and W. A. Kirk,
 {\it ``Topics in metric fixed point theory''},
 Cambridge Studies in Advanced Mathematics 28,
 Cambridge University Press (1990).

\bibitem{REF:Gossez_LamiDazo1972_Pacific}
 J.-P. Gossez and E. Lami Dazo,
 {\it ``Some geometric properties related to the fixed point theory for
 nonexpansive mappings''},
 Pacific J.\ Math., {\bf 40} (1972), 565--573.

\bibitem{REF:Hirano1982_JMSJapan}
 N. Hirano,
 {\it ``Nonlinear ergodic theorems and weak convergence theorems''},
 J.\ Math.\ Soc.\ Japan, {\bf 34} (1982), 35--46.

\bibitem{REF:Lim1974_Pacific}
 T. C. Lim,
 {\it ``A fixed point theorem for families of nonexpansive mappings''},
 Pacific J.\ Math., {\bf 53} (1974), 487--493.

\bibitem{REF:Miyadera_Kobayasi1982_NATMA}
 I. Miyadera and K. Kobayasi,
 {\it ``On the asymptotic behaviour of almost-orbits of nonlinear contraction
 semigroups in Banach spaces''},
 Nonlinear Anal., {\bf 6} (1982), 349--365.

\bibitem{REF:Opial1967_BullAMS}
 Z. Opial,
 {\it ``Weak convergence of the sequence of successive approximations for
 nonexpansive mappings''},
 Bull.\ Amer.\ Math.\ Soc., {\bf 73} (1967), 591--597.

\bibitem{REF:TS2003_ProcAMS}
 T. Suzuki,
 {\it ``On strong convergence to common fixed points of nonexpansive semigroups
 in Hilbert spaces''},
 Proc.\ Amer.\ Math.\ Soc., {\bf 131} (2003), 2133--2136.

\bibitem{REF:Suzuki-thm}
 T. Suzuki,
 {\it ``Some remarks on the set of common fixed points
 of one-parameter semigroups of nonexpansive mappings
 in Banach spaces with the Opial property''},
 submitted.

\bibitem{REF:TSP_mcbt_1_04}
 T. Suzuki and W. Takahashi,
 {\it ``Strong convergence theorems of Mann's type for one-parameter
 nonexpansive semigroups in general Banach spaces''},
 submitted.

\bibitem{REF:Takahashi_ybook}
 W. Takahashi,
 {\it ``Nonlinear Functional Analysis''},
 Yokohama Publishers,
 Yokohama, 2000.

\end{thebibliography}
\end{document}